\documentclass[11pt]{amsart}

\usepackage{amsfonts,amsmath,amssymb,amsthm,mathrsfs,url,color,latexsym,rawfonts,lscape,graphicx}
\usepackage{hyperref}

\theoremstyle{plain}
  \newtheorem{theorem}{Theorem}[section]
  \newtheorem{lemma}{Lemma}[section]
  \newtheorem{proposition}{Proposition}[section]

  \newtheorem{definition}{Definition}[section]
  \newtheorem{remark}{Remark}[section]
  \newtheorem{example}{Example}[section]

   \newcommand{\beqn}{\begin{eqnarray}}
   \newcommand{\eeqn}{\end{eqnarray}}
   \newcommand{\beqs}{\begin{eqnarray*}}
   \newcommand{\eeqs}{\end{eqnarray*}}
   \newcommand{\ban}{\begin{eqnarray*}}
   \newcommand{\nan}{\end{eqnarray*}}
   \newcommand{\beq}{\begin{equation}}
   \newcommand{\eeq}{\end{equation}}
   
\newcommand{\Td}{\mathbb T^d}
\newcommand{\dt}{\partial_t}
\newcommand{\demi}{\frac{1}{2}}
\let\cal=\mathcal

  \newcommand{\RR}{{\mathbb R}}

  \newcommand{\R}{\RR}

\newcommand{\p}{\partial}

\newcommand{\Om}{\Omega}

\renewcommand{\det}{\mbox{det}}

  \newcommand{\diam}{{\mbox{diam}}}

  \numberwithin{equation}{section}
  \numberwithin{figure}{section}

\newcommand{\Rd}{\mathbb R^d}
\newcommand{\bm}{\mathbf{m}}

\begin{document}

\title[Optimal transport with discrete long range mean field interactions]
{Optimal transport with discrete long range mean field interactions}

%\date\today

\author[J. Liu]
{Jiakun Liu}
\address
	{School of Mathematics and Applied Statistics,
	University of Wollongong,
	Wollongong, NSW 2522, AUSTRALIA}
\email{jiakunl@uow.edu.au}

\author[G. Loeper]
{Gr\'egoire Loeper}
\address
{School of Mathematical Sciences, 
Monash University,
Melbourne, VIC 3800, AUSTRALIA}

\email{gregoire.loeper@monash.edu}

%\thanks{This work is supported by ARC DP170100929 ...}

%\subjclass[2000]{35J60, 35B45; ????, ????}

%\keywords{Optimal transport, Mean field interaction, Monge-Amp\`ere equation}

\begin{abstract}
We study an optimal transport problem where, at some intermediate time, the mass is accelerated by either an external force field, or self-interacting. We obtain regularity of the velocity potential, intermediate density, and optimal transport map, under conditions on the interaction potential that are related to the so-called Ma-Trudinger-Wang condition from optimal transport \cite{MTW}. 
\end{abstract}

\maketitle

\baselineskip=16.4pt
\parskip=3pt

\section{Introduction}

\subsection{Motivations}

The optimal transport problem goes back to a cost-minimisation problem in civil engineering proposed by Monge \cite{Monge},  later generalised to a class of optimisation problems by Kantorovich \cite{kanto1, kanto2}, with an elegant economic interpretation.
Later, through the contribution of Brenier \cite{Bre},  Benamou and Brenier \cite{BB}, Frisch and co-authors \cite{BF, FM, MTF}, the second author \cite{Lo06}, and Lee and McCann \cite{LM}, the connection between optimal transport and classical mechanics has also appeared very naturally. 
Indeed, a natural approach in Hamiltonian mechanics is to look for critical points of the action of a Lagrangian in order to find the evolution equation of the system. In some cases where the Lagrangian has some form of coercivity (as in the natural action $\int_{[0,T]\times \R^d} \rho(t,x)|v|^2(t,x)dtdx$, cf. \cite{BB}), critical points of the action can be obtained by minimisation, which is a natural formulation as an optimal transport problem. The continuous time formulation of the problem in terms of  curves on a space of probability measures is extensively addressed in the book by Ambrosio, Gigli and Savar\'e \cite{AGS}, and also in the books by Villani \cite{Vi03,Vi09}.

One of the main differences between the {\it economic} point of view and the {\it mechanical} point of view, is the addition of the time variable. 
The original transport problem starts with a cost function $c(x,y)$ that depends only on the starting and arrival points. The action minimising problem looks for curves or vector fields (depending on whether one uses the Lagrangian or Eulerian point of view). We shall speak either of the {\it point to point} or  {\it time continuous} problem to distinguish between the two situations.

 Although in both cases  the existence of optimisers rests now on a well established theory, the question of their regularity is relatively well understood for the point to point problems, while it is still largely unexplored in the time continuous case (see references below). In the point to point problem, by regularity we mean the smoothness of the optimal map sending one distribution of mass to the other; in the time continuous case, such a map usually exists too, but is more difficult to characterise. The main reason for that difference is that the regularity of the point to point problem relies on the study of an associated Monge-Amp\`ere equation, which is not always accessible in the time continuous case. 

\subsection{An example: The reconstruction problem in cosmology}
As an illustrative example, we go back to a previous work by the second author \cite{Lo06}, where he studies the motion of self-gravitating matter,  classically described by the Euler-Poisson system, for an application in cosmology, known as the {\it reconstruction problem}, a problem that has received a lot of attention in cosmology (see also \cite{BF, FM, MTF} and the references therein). Let us recall the model here:
a continuous distribution of matter with density $\rho$, moves along a velocity field $v$, and is accelerated by a gravitational field, itself given as the gradient of a potential $p$, linked to $\rho$ through the Poisson equation. 
The system is thus the following 
	\begin{equation}\label{EPs}
	\left\{\begin{array}{rl}
		\partial_t\rho + \nabla\cdot(\rho v) \!&=\ 0, \\
		\partial_t(\rho v) + \nabla\cdot(\rho v\otimes v) \!&=\ -\rho\nabla p, \\
		\Delta p \!&=\ \rho.
	\end{array}
	\right.
	\end{equation}
The first equation is the conservation of mass,  the second equation states that the acceleration field is given by $-\nabla p$, and the third equation is the Poisson coupling between gravitation and matter.

The reconstruction problem is to find a solution to \eqref{EPs} satisfying
\begin{equation*}
		 \rho|_{t=0} = \rho_0, \qquad \rho|_{t=T} = \rho_T,
	\end{equation*}
i.e. given initial and final densities, as opposed to the Cauchy (or initial value) problem, where one is given initial density and velocity.	
In \cite{Lo06} (see also \cite{Bre} for the case of the incompressible Euler equations), the reconstruction problem was formulated into a minimisation problem, minimising the action of the Lagrangian which is a convex functional of properly chosen variables. Through this variational formulation, the reconstruction problem becomes very similar to the time continuous formulation of the optimal transport problem of Benamou and Brenier \cite{BB}. 
Moreover, through the study of a dual problem, reminiscent of Monge-Kantorovich duality, partial regularity results for the velocity and the density were obtained.

The optimal transport problem of \cite{Lo06} was formulated as finding minimisers of the action
\beqn\label{defI}
I(\rho,v,p)=\demi\int_{0}^{T}\int_{\Td}  \rho(t,x)|v(t,x)|^2 + |\nabla p(t,x)|^2\, dx dt,
\eeqn
over all $\rho, p, v$ satisfying
\beqs
\dt \rho + \nabla\cdot(\rho v)=0,\\
\rho(0)=\rho_0,\quad  \rho(T)=\rho_T,\\
\Delta p = \rho-1,
\eeqs
where $\Td$ denotes the $d$-dimensional torus, as the study in \cite{Lo06} was performed in the space-periodic case.

\subsection{Goal of the paper} 

In this paper we address the problem of 
finding minimisers for the action
\beqn\label{defIgen}
I(\rho,v,p)=\int_{0}^{T}\left(\int_{\Rd}  \demi\rho(t,x)|v(t,x)|^2 dx + {\cal F}(t,\rho(t))\right) dt,
\eeqn
for more general ${\cal F}$, and we will obtain some partial results in that direction. 

Apart from the application in cosmology, several authors have looked at continuous optimal transport with or without interaction, notably through their natural connections with mean-field games. We refer the reader to the chapter 
\cite{BCS} and the references therein, where the connection is explored. Again, we also refer to the book \cite{AGS} where the notion of gradient flows on space of measures and its relation to optimal transport is addressed.

The present work is about the study of regularity of minimisers to \eqref{defIgen}. In \cite{LM}, Lee and McCann address the case where
\beqs
{\cal F}(t,\rho)=\int \rho(t,x) Q(t,x) dx,
\eeqs
which is obviously linear in $\rho$.
This Lagrangian corresponds to the case of a continuum of matter evolving in an external acceleration field given by $\nabla Q(t,x)$. We call this the non-interacting case for obvious reasons. This can be recast as an optimal transport problem where the cost function is given by
\beqn
c(x,y)=\inf_{\gamma(0)=x, \gamma(T)=y} \int_0^T \demi |\dot \gamma(t)|^2 + Q(t,\gamma(t)) dt,
\eeqn
where $\gamma$ is a smooth curve connecting $x$ and $y$. 
By assuming that $Q(t,x)=\varepsilon V(x)$ for some $V$ satisfying the structure condition
\beqs
-\int_0^1\int_0^\tau \left.\langle u,(1-t)\partial_s^2\mbox{Hess}\, V_{x+t(v+sw)}u \rangle \right|_{s=0}dt d\tau \geq C
\eeqs
for a constant $C>0$, for all $(x,v)$ in the tangent bundle $T\Td$ and for all unit tangent vectors $u, w$ in the tangent space $T_x\Td$ that are orthogonal to each other, 
Lee and McCann obtain that for  a small enough $\varepsilon>0$, the cost $c$ satisfies the conditions found in \cite{MTW} to ensure the regularity of the optimal map. (Note that in order to be consistent with our notations, here changed the sign of the potential actually considered in \cite{LM}.)

\subsubsection{The non-interacting, discrete case, and the associated Monge-Amp\`ere equation}

In the present paper we will restrict ourselves to the case where the force field only acts at a single discrete time between $0$ and $T$: 
\beqs
Q(t,x)=\delta_{t=T/2}Q(x).
\eeqs
The minimisation problem therefore becomes
\beqn\label{defIdisc}
I(\rho,v,p)=\int_{0}^{T}\int_{\R^d} \demi \rho(t,x)|v(t,x)|^2 dx dt + \int_{\R^d}\rho(T/2,x)Q(x) dx,
\eeqn
for some potential $Q$.
This will allow to remove the smallness condition on $Q$. 
\subsubsection{The mean-field case} 
We will be able to extend our result to the non-linear situation where the force field is given by
\beqn
\nabla Q(x) = \int \rho(t,y)\nabla \kappa(x-y) dy,
\eeqn
still acting at a single intermediate time.
This corresponds to the case where a particle located at $y$ accelerates another particle located at $x$ with an acceleration equal to $\nabla \kappa(x-y)$.  In this case we will show that the action to minimise becomes 
\begin{equation}\label{gravit}
I(\rho,v) = \frac12\int_0^T\int_{\R^d}  \rho|v|^2 dxdt + \frac12 \int_{\R^d}\int_{\R^d} \rho(T/2,x)\kappa(x-y)\rho(T/2,y)dxdy.
\end{equation}

Reasoning formally, one sees straight away that on $[0,T/2]$ we are solving the usual optimal transport problem in its ``Benamou-Brenier" formulation \cite{BB}, as well as on $[T/2,T]$, and therefore particles will travel with constant velocity in those two intervals.
At $t=T/2$, the velocity $v$ will be discontinuous.
 We will give a sufficient condition on $\kappa$ to ensure a smooth transport map and intermediate density.
{Unfortunately, the gravitational case, which corresponds to the Coulomb kernel 
\beq\label{Coulomb}
\kappa(x-y)= \frac{c_d}{|x-y|^{d-2}}
\eeq does not satisfy our condition.}

\subsection{Statement of the problem}
We start by giving a three formulations of the problem that will turn out to be equivalent.
\vspace{10pt}

{\bf Problem 1}
We consider ${\cal P}_2(\Rd)$ (in short ${\cal P}_2$) the set of probability measures on $\R^d$ with finite second moment, and a functional ${\cal F}:{\cal P}_2\to \R\cup{+\infty}$. For a given curve $\rho(t):[0,T]\to {\cal P}_2$ and a vector field $v(t) \in L^2\big([0,T]\times \R^d,dt\times d\rho(t)\big)$ we consider the action 
\beqn\label{eq:min1}
I(\rho,v) = \int_{[0,T]\times \R^d} \frac12 d\rho(t,x)|v(t,x)|^2 dt + {\cal F}(\rho(T/2))
\eeqn
and the constraints
\beqn 
&&\dt\rho + \nabla\cdot(\rho v)=0,\label{eq:cont}  \\
&&\rho(0)=\rho_0, \;\rho(T)=\rho_T \label{eq:boundary}.
\eeqn 
The Problem 1 is to minimise $I$ among all $\rho,v$ satisfying \eqref{eq:cont}--\eqref{eq:boundary}.

\vspace{10pt}

{\bf Problem 2}
Consider the space of continuously differentiable curves $\Gamma = C^1([0,T]; \R^d)$. To each $x\in \R^d\cup {\text{support}(\rho_0)}$ we associate $\gamma(t,x)\in \Gamma$ such that $\gamma(0,x)=x$. We consider 
\beqn\label{eq:min2}
J(\gamma) = \int_0^T \demi |\dt \gamma(t,x)|^2 dt\, d\rho_0(x) + {\cal F}(\rho(T/2)),
\eeqn
where\footnote{The subscript $_\#$ classically means that $\gamma(T/2)$ pushes forward $\rho_0$ onto $\rho(T/2)$.}
\beqn\label{eq:contcurve}
\rho(T/2)=\gamma(T/2)_\# \rho_0,
\eeqn
and under the constraint 
\beqn\label{eq:gammaT}
\rho(T)=\gamma(T)_\# \rho_0.
\eeqn
The Problem 2 is to minimise $J$ among all $\gamma$ satisfying \eqref{eq:contcurve}--\eqref{eq:gammaT}.

\vspace{10pt}

{\bf Problem 3} 
For $\mu,\nu \in {\cal P}_2$, letting $\Gamma_{\mu,\nu}$ be the set of probability measures on $\R^d\times\R^d$ with marginals $\mu$ and $\nu$, we recall that the Wasserstein distance (of order 2) between $\mu$ and $\nu$ is given by
\beqn\label{def:W2}
W_2^2(\mu,\nu)=\inf_{\pi\in \Gamma_{\mu,\nu}} \int_{\R^d\times\R^d}\demi|x-y|^2 d\pi(x,y).
\eeqn
We define 
\beqn\label{defK}
K(\rho) = \frac2T W_2^2(\rho_0,\rho) + {\cal F}(\rho) + \frac2T W_2^2(\rho,\rho_T).
\eeqn
The Problem 3 is to minimise $K$ among all $\rho \in {\cal P}_2$. It co\"incides with the notion of  Wasserstein Barycenters, see \cite{AgCa}, \cite{KimPass}, when $${\cal F}(\rho) = CW_2^2(\rho_I,\rho),\quad C>0$$ for some intermediate measure $\rho_I\in {\cal P}_2$.

\vspace{10pt}

From classical results of optimal transport (see \cite{Vi03, Vi09, AGS}), in the case where ${\cal F(\rho)}$ is convex and lower semi continuous (l.s.c.) in $\rho$ there holds:
\begin{proposition}\label{prop:equiv}
Let $\rho_0,\rho_T \in {\cal P}_2\cap L^1(\Rd)$. 
Assume that ${\cal F(\rho)}$ is convex and lower semi continuous (l.s.c.) in $\rho$, and that Problem 1,2 or 3 has at least one admissible solution, then Problems 1,2,3 are equivalent, and moreover there holds
\beqn
v(0,x)=\dt\gamma(0,x)=\frac{2}{T}(\nabla\Phi^*(x)-x)=:\nabla\phi(x),
\eeqn where $v,\gamma$ are respectively from Problems 1,2 and 
$\Phi^*$ is a convex potential such that $\nabla\Phi^*_\#\rho_0=\rho$ and $\rho$ is the optimiser in Problem 3.
\end{proposition}

The second proposition gives the Euler-Lagrange equation characterising the optimiser. It is based on the Riemannian metric induced on ${\cal P}_2$ by $W_2$ (see again the above references for a complete coverage).

\begin{proposition}
Let $\rho$ be the optimiser in Problem 3. There exists a vector field $w=\nabla\Xi \in L^2(d\rho)$ and two convex potentials $\Phi, \Psi$ such that\footnote{Note that $\Phi^*$ above is the Legendre transform of $\Phi$.}
\beqn
\nabla\Phi_\# \rho \!\!&=&\!\! \rho_0,\\  
\nabla\Psi_\# \rho \!\!&=&\!\! \rho_T,\\
\nabla\Psi+\nabla\Phi \!\!&=&\!\! 2x+ \frac{T}{2}\nabla\Xi,
\eeqn
moreover $\nabla\Xi$ can be identified to the gradient of ${\cal F}$ with respect to the Wasserstein metric, i.e. for all curve $\rho_t \subset {\cal P}_2\cup Dom({\cal F})$ passing through $\rho$ at $t=0$, and such that 
\ban
\dt \rho_t + \nabla\cdot(\rho_t v_t)=0,
\nan
for some $v_t \in L^\infty([0,T]; L^2(d\rho_t))$
there holds $$\frac{d}{dt}{\cal F}(\rho_t)\Big|_{t=0} = \int d\rho\, \nabla\Xi \cdot v_0.$$ 
We will denote $\nabla\Xi=\nabla_W{\cal F}\Big|_{\rho}$.
\end{proposition} 

We next characterise $\Xi$ in the two cases of interest for the present paper:
\begin{proposition}
\begin{itemize}
\item[i)] Assume that ${\cal F}(\rho)$ is as in \eqref{defIdisc}. Then $$\Xi= \frac{Q}{T}.$$
\item[ii)] If ${\cal F}(\rho)$ as in \eqref{gravit} then
$$\Xi = \frac{1}{T}\int_{\Rd} \kappa(\cdot-y) d\rho(y)=:\frac{Q_\rho}{T}.$$
\end{itemize}
We can also completely characterise the optimal velocity: with $\phi$ as in Proposition \ref{prop:equiv}, 
\begin{itemize}
\item[$(i)$] when $t\in[0,\frac{T}{2})$, $v(t,x)=\nabla\phi(x_0)$, where $x=x_0+t\nabla\phi(x_0)$; $\rho(t)=(x+t\nabla\phi)_\# \rho_0$.
\item[$(ii)$] when $t\in(\frac{T}{2}, T]$, $v(t,x)=\nabla\phi(x_0) + \nabla Q(z)$, where $x=z + (t-\frac{T}{2})(\nabla\phi(x_0) + \nabla Q(z))$ and $z=x_0+\frac{T}{2}\nabla\phi(x_0)$; $\rho(t)=\big(x+(t-T/2)(\nabla\phi+ \nabla Q)\big)_\# \rho_{T/2}$.
\end{itemize}
Finally the optimal map: $\mathbf{m}(x)$ such that
$\mathbf{m}(\gamma(0,x)) = \gamma(T,x)$ (for $\gamma$ the optimiser in Problem 2)  will be given by
\beqn\label{def:bm}
\mathbf{m}(x) = x + T\nabla\phi + \frac{T}{2}\nabla Q \left(x+\frac{T}{2}\nabla\phi \right),
\eeqn 
and there also holds 
\beqn\label{def:bm2}
\mathbf{m}(x)=\nabla\Psi(\nabla\Phi^*).
\eeqn
\end{proposition}

\subsection{Assumptions}

From the previous observations, and in order to motivate our assumptions, let us derive formally the equation giving the initial velocity potential $\phi$.
Let the initial density $\rho_0$ be supported on a bounded domain $\Omega_0\subset\mathbb{R}^d$, and the final density $\rho_T$ be supported on a bounded domain $\Omega_T\subset\mathbb{R}^d$, satisfying the balance condition
	\begin{equation}\label{bala}
		M:=\int_{\Omega_0} \rho_0(x)\,dx = \int_{\Omega_T} \rho_T(y)\,dy.
	\end{equation}

Starting from the definition \eqref{def:bm}, we introduce the modified potential functions	
	\begin{equation}\label{modphi}
		\tilde\phi(x):=\frac{T}{2}\phi(x)+\frac12|x|^2,\qquad\mbox{ and  }\quad \tilde Q(z) := \frac{T}{2}Q(z)+|z|^2.
	\end{equation}	
By computing the determinant of the Jacobian $D\mathbf{m}$ and noting that $\mathbf{m}_\#\rho_0=\rho_T$, i.e. that $\mathbf{m}$ pushes forward the measure $\rho_0$ onto the measure $\rho_T$, one can derive the Monge-Amp\`ere type equation (see \S\ref{s2} for detailed computation)
	\begin{equation}\label{MA0}
		\det\,\left[ D^2\tilde\phi - \left(D^2\tilde Q(\nabla\tilde\phi)\right)^{-1}\right] = \left(\frac{1}{\det D^2\tilde Q{(\nabla\tilde\phi)}}\right) \frac{\rho_0}{\rho_T\circ\mathbf{m}},
	\end{equation}
with a natural boundary condition 
	\begin{equation}\label{bc00}
		\mathbf{m}(\Omega_0)=\Omega_T.
	\end{equation}	 

To ensure the regularity of the solution $\tilde\phi$ (equivalently that of $\phi$) to the boundary value problem \eqref{MA0} and \eqref{bc00}, 
it is necessary to impose certain conditions on the {potential} function $\tilde Q$ (equivalently on $Q$) and the domains $\Om_0, \Om_T$.  
In this paper we assume $\tilde Q$ satisfies the following conditions: 
\begin{itemize}
\item[\textbf{(H0)}] The {potential} function $\tilde Q$ belongs to  $C^4(\Rd)$.

\item[\textbf{(H1)}] The {potential} function $\tilde Q$ is uniformly strictly convex, namely $D^2\tilde Q \geq \varepsilon_0 I$ for some $\varepsilon_0 >0$.
\item[\textbf{(H2)}] The {potential} function $\tilde Q$ satisfies for all $\xi, \eta\in\mathbb{R}^d$ with $ \xi \bot \eta$, 
	\begin{equation}\label{newH20}
		\sum_{i,j,k,l,p,q,r,s}\left(\tilde Q_{ijrs} - 2\tilde Q^{pq}\tilde Q_{ijp} \tilde Q_{qrs} \right)\tilde Q^{rk}\tilde Q^{sl} \xi_k\xi_l\eta_i\eta_j \leq -\delta_0|\xi|^2|\eta|^2,
	\end{equation}
where $\{\tilde Q^{ij}\}$ is the inverse of $\{\tilde Q_{ij}\}$, and $\delta_0$ is a positive constant.
When $\delta_0=0$, we call it \textbf{(H2w)}, a weak version of {(H2)}.
\end{itemize}
Note that conditions (H0) and (H1) imply that the inverse matrix $(D^2\tilde Q)^{-1}$ exists, and ensure that equation \eqref{MA0} is well defined. Condition (H2) is an analogue of the Ma-Trudinger-Wang (MTW) condition in optimal transportation, which is necessary for regularity results (note the factor 2 however). 
We shall give more detailed interpretations and examples in \S\ref{s4}.

\vspace{10pt}

\section{Results}
Our first main result is the following
\begin{theorem}\label{thm:equiv}
Under assumptions (H0) and (H1) Problem 1,2,3 are equivalent to solving an optimal transport problem with cost function $c(x,y)=\tilde Q^*(x+y)$, where $\tilde Q^*$ is the Legendre transform of $\tilde Q$ given in \eqref{modphi}. There exists a potential $\tilde \phi$ such that the optimal transport map of $\rho_0$ onto $\rho_T$ with cost $c$ will be given by $T(x) = y$ such that $D_xQ^*(x+y)=D_x\tilde\phi(x)$. Then $\phi=\frac2T(\tilde\phi-|x|^2/2)$ will be the initial velocity potential as in Proposition \ref{prop:equiv}.
\end{theorem}

Our next result is a regularity result:
\begin{theorem}\label{mainthm}
Let $\phi$ be the {initial} velocity potential as above.
Assume the {potential} function $\tilde Q$ satisfies conditions (H0), (H1) and (H2), $\Omega_T$ is $q$-convex with respect to $\Omega_0$ (defined in \S\ref{s3}).
Assume that $\rho_T\geq c_0$ on $\Omega_T$ for some positive constant $c_0$, $\rho_0\in L^p(\Omega_0)$ for some $p>\frac{d+1}{2}$, and the balance condition \eqref{bala} is satisfied.
Then, the velocity potential $\phi$ is $C^{1,\alpha}(\overline\Omega_0)$ for some $\alpha\in(0,1)$.

If furthermore, $\Om_0, \Om_T$ are $C^4$ smooth and uniformly $q$-convex with respect to each other, $\rho_0\in C^2(\overline\Om_0), \rho_T\in C^2(\overline\Om_T)$, then $\phi\in C^3(\overline\Om_0)$, and higher regularity follows from the theory of linear elliptic equations. In particular, if $\tilde Q, \Om_0, \Om_T, \rho_0, \rho_T$ are $C^\infty$, then the velocity potential $\phi$ is in $C^\infty(\overline\Om_0)$.
\end{theorem}

The proof of Theorem \ref{mainthm} follows from Theorem \ref{thm:equiv}, and from the observation that condition (H2) is equivalent to $\tilde Q^*$ satisfying the MTW condition.
Under this formulation, the regularity results then follow from the results in \cite{Lo06, Liu, MTW, TW1}. See also \cite{CW, FKM} for $C^{1,\alpha}$ regularity results under the condition (H2w) and some additional conditions on domains.

As a by-product of those two results we obtain the following:
\begin{theorem}
Consider the optimal transport problem with cost $c(x,y)=R(x+y)$ for some $R:\R^d\to \R$ convex. Then this problem is equivalent to the minimisation problem \eqref{defIdisc} with potential $Q(z) = \frac2T (R^*(z)-|z|^2)$, for $R^*$ the Legendre Fenchel transform of $R$.
\end{theorem}

%%%%%%%%%%%%%%%%%%%%%%%%%

For the mean-field case, we have the following existence and uniqueness result:
\begin{theorem}\label{theo:2.5}
Assume that $\kappa\in C^4(\R^d;\R^+)$, is convex.
There exists a unique minimiser to 
problem \eqref{gravit}. Moreover, once $\rho(T/2)$ is known, letting $\bar \rho=\rho(T/2)$, the minimiser will be the same as the solution of non-interacting problem \eqref{defIdisc} where $Q$ is given by
\beqn\label{newpo}
Q(x) = \int_{\R^d} \bar\rho(y) \kappa(x-y) dy.
\eeqn
\end{theorem}

Under additional assumptions on the kernel $\kappa$ and the domains, we have the following regularity result:
\begin{theorem}\label{mainthm2}
Let $\Omega_{T/2}:=\nabla\tilde\phi(\Om_0)$ with $\tilde\phi$ given in \eqref{modphi}.
Assume moreover that $\kappa$ satisfies 
\begin{itemize}
\item[(\bf{H2c})] for any $\xi, \eta\in\mathbb{R}^d$, $x,y\in\Omega_{T/2}$,
	\begin{equation*}
		\sum_{i,j,k,l,p,q,r,s}\left(D^4_{ijrs}\kappa(x-y)\right) \tilde\kappa^{rk}\tilde\kappa^{sl}\xi_k\xi_l\eta_i\eta_j \leq 0,
	\end{equation*}
where $\{\tilde\kappa^{ij}\}$ is the inverse of $\{\kappa_{ij}+\frac{2}{T}I\}$.
\end{itemize}
Assume that $\Omega_0$, $\Omega_T$ are smooth convex domains. Let $\varphi\in C^2(\overline\Om_0)$, $\psi\in C^2(\overline\Om_T)$ be convex defining functions of $\Om_0$, $\Om_T$, respectively. Suppose that for any $z, w\in\Om_{T/2}$,
\beq\label{cx1}
\left[\varphi_{ij}(x)+\frac{TM}{8}\kappa_{ijk}(z-w)\varphi_k(x)\right] \geq b_0\delta_{ij} \qquad \forall x\in\p\Om_0, 
\eeq
and
\beq\label{cx2}
\left[\psi_{ij}(y)+\frac{TM}{8}\kappa_{ijk}(z-w)\psi_k(y)\right] \geq b_1\delta_{ij} \qquad \forall y\in\p\Om_T,
\eeq
where $b_0, b_1$ are two constants, and $M$ is the total mass in \eqref{bala}. 

If $b_1\geq0$, $\Om_T$ is $q$-convex with respect to $\Om_0$, where $Q$ is given by \eqref{newpo}, thus the first conclusion of Theorem \ref{mainthm} holds. 
If furthermore, $b_0, b_1 > 0$ are positive, $\Om_0, \Om_T$ are uniformly $q$-convex with respect to each other, and thus the second conclusion of Theorem \ref{mainthm} holds, namely the initial velocity potential $\phi$ is smooth provided $\kappa, \Om_0, \Om_T, \rho_0, \rho_T$ are smooth, which in turn implies that the intermediate density $\rho(T/2)$ is smooth. 
\end{theorem}

The proof of Theorem \ref{mainthm2} relies on the observation that the $q$-convexity and the condition (H2c) is preserved under convex combinations, and therefore by convolution with the density $\rho(T/2)$, and on some a priori $C^1$ estimates on the potential.

\vspace{5pt} 
The paper is organised as follows: 
In \S\ref{s2}, we derive equation \eqref{MA0} formally by straightforward computations. 
In \S\ref{s3}, we introduce the two-step optimal transport problem and prove Theorem \ref{thm:equiv}. By assuming conditions (H0)--(H1) we have the existence and uniqueness of the velocity potential $\phi$. Moreover, we provide an interpretation of the cost function from a natural mechanical point of view. 
In \S\ref{s4}, we introduce the condition (H2), which is crucial in obtaining the regularity of $\phi$. 
In \S\ref{s5}, upon formulating our reconstruction problem into an optimal transport problem, we have the regularity of $\phi$ and conclude Theorem \ref{mainthm}.
In \S\ref{s6}, we consider the mean-field case under the condition (H2c), which is preserved by convex combinations, and then prove Theorem \ref{theo:2.5} and \ref{mainthm2}.

%%%%%%%%%%%%%%%%%%%%%%%%%%%%%%%%%%%%%%%%
%%%%%%%%%%%%%%%%%%%%%%%%%%%%%%%%%%%%%%%%
\vspace{10pt}
\section{Formal derivation of equation \eqref{MA0}}\label{s2} 
Throughout the following context, unless mentioned otherwise, the function $\phi$ always denotes the initial velocity potential, namely at time $t=0$, the velocity
	\begin{equation}\label{realpo}
		v(0,x) = \nabla \phi(x),\quad \mbox{for } x\in\Omega_0.
	\end{equation}
In order to derive the equation for $\phi$, let's track a single point $x\in\Omega_0$. 

Recalling that at $t=T/2$, the {potential} $\nabla Q$ affects the velocity $v=\nabla\phi$, the final point $y=\mathbf{m}(x)$ is given by
	\begin{equation}\label{2step}
		\mathbf{m}(x) = x + T\nabla\phi + \frac{T}{2}\nabla Q \left(x+\frac{T}{2}\nabla\phi \right).
	\end{equation}
The Jacobian matrix of $\bm$ is
	\begin{equation*}
	\begin{split}
		D\mathbf{m} &= I + T D^2\phi + \frac{T}{2} D^2Q \cdot \left[ I + \frac{T}{2}D^2\phi \right] \\
				&= \left[ I + \frac{T}{2}D^2\phi \right] \cdot \left[ I + \frac{T}{2}D^2Q \right] + \frac{T}{2}D^2\phi,
	\end{split}
	\end{equation*}
where $I$ is the $d\times d$ identity matrix, {and the matrix $(D^2Q)$ is taken at $(x+T/2\nabla\phi$)}.

Define	
	\begin{equation*}
		A := \left[ I + \frac{T}{2}D^2Q \right],
	\end{equation*}
and assume that the matrix $(I+A)$ is invertible. Then
	\begin{equation}\label{Jac}
	\begin{split}
		D\mathbf{m} &= \left[ I + \frac{T}{2}D^2\phi \right] \cdot A + \frac{T}{2}D^2\phi \\
				&= \frac{T}{2}D^2\phi \cdot \left[ I + A\right] + A \\
				&= \frac{T}{2}D^2\phi \cdot \left[ I + A\right] + (I+A) - I \\
				&= [I+A]\left[ \frac{T}{2}D^2\phi + I - (I+A)^{-1} \right].
	\end{split}
	\end{equation}	
Computing the determinant of $D\mathbf{m}$, we have
	\begin{equation*}
	\begin{split}
		\det\,D\mathbf{m} &= \det\, \left[ I + A\right]\, \det\, \left[\frac{T}{2}D^2\phi + I - (I+A)^{-1} \right] \\
				&= \det\, \left[ I + A\right]\, \det\, \left[\frac{T}{2}D^2\phi + I - \left( 2I + \frac{T}{2}D^2Q\big(x+ \frac{T}{2}\nabla\phi(x)\big) \right)^{-1}\right]. 
	\end{split}
	\end{equation*}
	
Recall that the modified velocity potential $\tilde\phi$ and {potential} function $\tilde Q$ are given by
	\begin{equation*}
		\tilde\phi(x)=\frac{T}{2}\phi(x)+\frac12|x|^2 \quad \mbox{ and }\quad \tilde Q(z) = \frac{T}{2}Q(z)+|z|^2,
	\end{equation*}	
we have $D^2\tilde Q = (I+A)$ and
	\begin{equation*}
		\det\,D\mathbf{m} = \det D^2\tilde Q \, \det \left[ D^2\tilde\phi - \left(D^2\tilde Q(\nabla\tilde\phi)\right)^{-1}\right].
	\end{equation*}
Therefore, we obtain the Monge-Amp\`ere equation
	\begin{equation*}
		\det\,\left[ D^2\tilde\phi - \left(D^2\tilde Q(\nabla\tilde\phi)\right)^{-1}\right] = \frac{\det\,D\mathbf{m}}{\det D^2\tilde Q}. 
	\end{equation*}
Note that $\mathbf{m}_\#\rho_0=\rho_T$ (defined in \eqref{pushfd}), thus $|\det\,D\mathbf{m}(x)|=\rho_0(x)/\rho_T(\mathbf{m}(x))$.
Then we obtain the equation \eqref{MA0},
	\begin{equation}\label{MA}
		\det\,\left[ D^2\tilde\phi - \left(D^2\tilde Q(\nabla\tilde\phi)\right)^{-1}\right] = \left(\frac{1}{\det D^2\tilde Q}\right) \frac{\rho_0}{\rho_T\circ\mathbf{m}},
	\end{equation}
with an associated natural boundary condition 
	\begin{equation}\label{bcon}
		\mathbf{m}(\Omega_0)=\Omega_T.
	\end{equation}

\begin{remark} 
\emph{In the continuous case \eqref{EPs}, the second author obtained in \cite{Lo06} partial regularity for $\phi$, that hold only in the interior (with respect to time) of the domain. In particular, there was no result regarding the initial velocity.
In this paper, we obtain the regularity of $\phi$ over the whole domain in the discrete case by using the regularity in optimal transportation. }
\end{remark}

We need to introduce a notion of weak solutions for equation \eqref{MA}.
\begin{definition}\label{Bresol}
A function $\tilde \phi$ is said to be a weak {\it Brenier}  solution to \eqref{MA} whenever $\mathbf{m}$ defined from $\tilde \phi$ as in \eqref{2step} is such that 
\begin{equation}\label{pushfd}
\mathbf{m}_\#\rho_0=\rho_T,
\end{equation}
namely, for all $B\subset \R^d$ Borel, there holds $\rho_0(\mathbf{m}^{-1}(B)) = \rho_T(B)$ (this is also called that $\mathbf{m}$ pushes forward $\rho_0$ onto $\rho_T$).
\end{definition} 

\begin{remark}
\emph{It is known that, depending on the geometry of the support of $\rho_T$, the notion of {\it Brenier} solution might be equivalent to the notion of {\it Aleksandrov} solution.  }
\end{remark}

%%%%%%%%%%%%%%%%%%%%%%%%%%%%%%%%%%%%%%%%
\vspace{10pt}
\section{two-steps transport}\label{s3}

Problem (\ref{defIdisc}) falls into the class of optimal transport problems with a general cost function
\begin{equation}\label{costLM}
c_T(x,y)=\inf\int_0^T L(\gamma(t), \dot{\gamma}(t))\,dt,
\end{equation}
where the infimum is taken over all smooth curves $\gamma(\cdot)$ satisfying $\gamma(0)=x$ and $\gamma(T)=y$,
as considered in \cite{LM}.
In our case the Lagrangian is defined by
 $L(x,v)=\frac12|v|^2+\delta_{t=T/2}Q(x)$.
Moreover, we can compute explicitly the optimal path $\gamma$ by dividing the transport map $\mathbf{m}=\mathbf{m_2}\circ\mathbf{m_1}$ such that at $t=\frac{T}{2}$,
	\begin{align}\label{map1}
		z &= \mathbf{m_1}(x) \\
			&= x+\frac{T}{2}\nabla\phi(x) = \nabla\tilde\phi(x), \nonumber
	\end{align}
and at $t=T$, 
	\begin{align}\label{map2}
		y &= \mathbf{m_2}(z) \nonumber \\
			&= z+ \frac{T}{2}\nabla Q(z) + \frac{T}{2}\nabla\phi(x) = 2z + \frac{T}{2}\nabla Q(z) - x \\
			&= \nabla\tilde Q(z) -x. \nonumber
	\end{align}
Correspondingly, 
	\begin{equation}\label{newcost}
	\begin{split}
		c_T(x,y) &= \inf_{z} \left\{ \frac{1}{2}\left|\frac{z-x}{T/2}\right|^2\cdot\frac{T}{2} + \frac{1}{2}\left|\frac{y-z}{T/2}\right|^2\cdot\frac{T}{2} + Q(z) \right\} \\
			&= \inf_{z} \left\{ \frac{1}{T}\left(|z-x|^2 + |y-z|^2\right) + Q(z) \right\}.
	\end{split}
	\end{equation}
Formally at the minimizer by taking $0=\frac{\partial c_T(x,y)}{\partial z}$ one can recover the equality \eqref{map2}, namely
	\begin{equation*}
		\frac{y-z}{T/2} = \frac{z-x}{T/2} + \nabla Q(z).
	\end{equation*}
Furthermore, through a straightforward computation, one has 	
	\begin{equation*}
	\begin{split}
		c_T(x,y) &= \inf_z\left\{ \frac{2}{T}|z|^2 - \frac{2}{T}z\cdot(x+y) + Q(z) \right\} + \frac{1}{T}\left(|x|^2+|y|^2\right) \\
			&= \frac{2}{T} \inf_z\left\{ |z|^2 - z\cdot (x+y) + \frac{T}{2}Q(z) \right\} + \frac{1}{T}\left(|x|^2+|y|^2\right).
	\end{split}
	\end{equation*}
	
Now, let $\tilde\phi^*, \tilde Q^*$ be Legendre transforms of convex functions $\tilde\phi, \tilde Q$ respectively, that is
	\begin{equation*}
	\left\{\begin{array}{ll}
		 (x, \nabla\tilde\phi(x)) &\sim\ (\nabla\tilde\phi^*(z), z), \\
		 (z, \nabla\tilde Q(z)) &\sim\ (\nabla\tilde Q^*(p), p),
	\end{array}
	\right.
	\end{equation*}
where,
	\begin{equation}\label{varp}
		p := \nabla\tilde Q(z) = y+x. 
	\end{equation}	 
Recall that $\tilde Q(z)=\frac{T}{2}Q(z)+|z|^2$ defined in \eqref{modphi}, one has
	\begin{equation*}
	\begin{split}
		c_T(x,y) &= \frac{2}{T} \inf_z \left\{ \tilde Q(z) - z\cdot p \right\} + \frac{1}{T}\left(|x|^2+|y|^2\right)\\
			&= -\frac{2}{T} \sup_z \left\{ z\cdot p - \tilde Q(z) \right\} + \frac{1}{T}\left(|x|^2+|y|^2\right)\\
			&= -\frac{2}{T} \tilde Q^*(p) + \frac{1}{T}\left(|x|^2+|y|^2\right),
	\end{split}
	\end{equation*}
where $\tilde Q^*$ is exactly the Legendre transform of $\tilde Q$ as defined above \eqref{varp}. 
{Note that $\tilde Q^*$ is well defined and $C^4$ smooth under assumptions (H0) and (H1).}

Note that the terms involving only $x$ or $y$ do no effect the  optimal transport map, therefore, we will look at an optimal transport problem with cost 
\begin{equation}\label{cost}
c(x,y)= \tilde Q^*(x+y),
\end{equation}
and seek to maximise the cost functional
	\begin{equation*}
		\mathcal{C}(\mathbf{s}) = \int_{\Omega_0} \rho_0(x)\tilde Q^*(x+\mathbf{s}(x))\,dx,
	\end{equation*}
among all maps $\mathbf{s} : \Omega_0 \to \Omega_T$ such that $\mathbf{s}_\#\rho_0=\rho_T$.
Note that here we consider the maximisation instead of the minimisation problem as in \cite{Lo09, MTW}.

When the cost function $c$ is strictly convex as is the case for $\tilde Q^*$ in \eqref{cost} satisfying (H0)--(H1), it was proved \cite{Caff, GM} that a unique optimal mapping can be determined by the potential functions, that leads to the Monge-Amp\`ere equation 
\begin{equation}\label{MA**} 
\det\,\left[ D^2\tilde\phi(x) - D^2\tilde Q^*(p) \right] = \left(\det D^2\tilde Q^*(p)\right) \frac{\rho_0(x)}{\rho_T\circ\mathbf{m}(x)},
\end{equation}
with the natural boundary condition \eqref{bcon}, where $p=x+y$. 
Note that  the matrix $(D^2\tilde\phi-D^2\tilde Q^*)$ is nonnegative, and $D^2\tilde Q^*$ is positive definite by condition (H1), which makes equation \eqref{MA**} elliptic. Note also that in the absence of regularity, one has to understand the solution to \eqref{MA**} in the weak ``Brenier" sense in Definition \ref{Bresol} (see \cite{Lo09}).

In our case, 
\beqn\label{Q*tildephi}
\nabla \tilde Q^*(x+y)=\nabla\tilde\phi(x),
\eeqn 
and using the properties of the Legendre transform
\begin{align*}
y&=\nabla\tilde Q(\nabla\tilde\phi(x))-x\\
&=\frac{T}{2}{\nabla}Q(\nabla\tilde\phi(x))+2\nabla\tilde\phi(x) - x.
\end{align*}
From the above computations, the initial velocity is given by $\frac{2}{T}(z-x)$ where $z=\nabla \tilde{Q}^*(x+y)$, hence following (\ref{Q*tildephi}) we recover that 
\beqn
v(0,x)=\nabla\phi(x) = \frac{2}{T}(\nabla\tilde{\phi}(x)-x),
\eeqn
while $\nabla\tilde{\phi}(x)=z$, where $z$ is the mid-point of the trajectory at $t=T/2$.

\vspace{5pt}
	
\begin{lemma}[$C^1$-bound]\label{C1lem}
Let $\Om_0, \Om_T$ be two bounded domains and $\tilde\phi$ be a solution of \eqref{MA**} and \eqref{bcon}.
Assume that the potential $\tilde Q$ satisfies {(H0)} and {(H1)}.
Then
	\begin{equation}\label{C1bdd}  
		|\nabla\tilde\phi(x)| \leq C,
	\end{equation}
for all $x\in\overline\Om_0$, where the constant $C$ depends on $\Om_0, \Om_T$ and $\tilde Q$. 
Furthermore, at time $t=T/2$, the intermediate domain $\Om_{T/2}:=\mathbf{m_1}(\Om_0)=\nabla\tilde\phi(\Om_0)$ is bounded.
\end{lemma}

\begin{proof}
From \eqref{Q*tildephi} and $p=x+y$, 
	\[ \nabla\tilde\phi(x) = \nabla\tilde Q^*(x+y), \]
where $\tilde Q^*=\tilde Q^*(p)$ is the Legendre transform of $\tilde Q=\tilde Q(z)$, thus is smooth and strictly convex. 
% $D^2\tilde Q \geq \varepsilon_0 I$ implies that $|D^2\tilde Q^*|<C$.
As $x\in\Om_0$, $y\in\Om_T$, $x+y\in B_R(0)$ for a bounded constant $R$.
Hence $|\nabla\tilde Q^*(x+y)|\le C$, and \eqref{C1bdd} is obtained. 

Recall that \eqref{map1}, at time $t=\frac{T}{2}$, $z=\mathbf{m_1}(x)=\nabla\tilde\phi(x)$. 
The inequality \eqref{C1bdd} implies that the intermediate domain $\Om_{T/2}=\mathbf{m_1}(\Om_0)$ is bounded. 
\end{proof}

\vspace{5pt}

We remark that from the uniqueness of $v$ in \cite{Lo06}, a solution of \eqref{MA**} is thus the velocity potential $\tilde\phi$. 
Therefore, we have  
\begin{proposition}\label{p31}
The two-steps gravitational transport problem is an optimal transportation associated with the cost function \eqref{cost}.
By assuming conditions (H0) and (H1), we have the existence and uniqueness (up to a constant) of solution $\tilde\phi$ to the boundary value problem \eqref{MA**} and \eqref{bcon} in the weak ``Brenier" sense. 

And in turn by \eqref{modphi}, we have the existence and uniqueness (up to a constant) of the initial velocity potential $\phi$.
\end{proposition}

Theorem \ref{thm:equiv} follows then directly from Proposition \ref{p31}.
Additionally, to prove Theorem \ref{mainthm} it is equivalent to obtain the regularity of solutions of \eqref{MA**} in optimal transportation, that requires appropriate convexity conditions on domains $\Omega_0, \Omega_T$, and more importantly some {structure} conditions on the {potential} function $\tilde Q$ to be described in the following sections.

\vspace{10pt}
\section{Conditions on the potential function}\label{s4}
 
From Proposition \ref{p31}, in order to obtain the regularity of the velocity potential $\tilde\phi$, it suffices to consider the optimal transportation with the cost function \eqref{cost}. 
From the results of \cite{Lo09}, it is now well understood that the so-called Ma-Trudinger-Wang (MTW) condition (introduced in \cite{MTW}) is necessary (at least in its weak form) for the regularity of optimal mappings.

Firstly, let us recall the MTW condition in optimal transportation.
For a general cost function $c(x,y) : \mathbb{R}^d\times\mathbb{R}^d\to\mathbb{R}$, use the notation
	\begin{equation*}
		c_{ij\cdots,kl\cdots} = \frac{\partial}{\partial x_i}\frac{\partial}{\partial x_j}\cdots\frac{\partial}{\partial y_k}\frac{\partial}{\partial y_l}\cdots c,
	\end{equation*}
and let $[c^{i,j}]$ denote the inverse of $[c_{i,j}]$.
The MTW condition is that 
	\begin{equation}\label{cMTW}
		MTW := \left( c_{ij,rs} - c^{p,q}c_{ij,p}c_{q,rs} \right) c^{r,k}c^{s,l}\xi_k\xi_l\eta_i\eta_j \geq c_0|\xi|^2|\eta|^2,
	\end{equation}
for any $\xi, \eta \in \mathbb{R}^d$ and $\xi\bot\eta$, where $c_0>0$ is a constant.
When $c_0=0$, it is called the weak MTW condition.

In our case, the cost function is given by the {potential} function $\tilde Q^*$ in \eqref{cost}. 
To introduce the analogous MTW condition, denote the matrix
	\begin{equation}\label{nota}
		D^2\tilde Q^*(p)  = \left(D^2\tilde Q(\nabla\tilde\phi)\right)^{-1} =: \mathcal{A}(z),
	\end{equation}
where $p=x+y$ and $z=\nabla\tilde\phi$, as in \eqref{varp} and \eqref{map1}, respectively.
Since $c(x,y)=\tilde Q^*(x+y)$, by differentiation we have
	\begin{equation}\label{dualMTW}
		MTW = \left( \tilde Q^*_{ijrs} - \tilde Q^{* pq}\tilde Q^*_{ijp}\tilde Q^*_{qrs} \right) \tilde Q^{*rk}\tilde Q^{*sl}\xi_k\xi_l\eta_i\eta_j.
	\end{equation}
From the Legendre transform and \eqref{nota}, one has
\beqs
	z=\nabla\tilde Q^*(p),\quad\mbox{ and }\quad \frac{\p z_i}{\p p_j}  = D^2_{ij}\tilde Q^*(p) = \mathcal{A}_{ij}(z).
\eeqs
Hence, using Einstein summation one has
	\begin{align}
		D^3_{ijr}\tilde Q^* &= \frac{dz_k}{dp_r} D_k\mathcal{A}_{ij} = \mathcal{A}_{kr} D_k\mathcal{A}_{ij},  \label{mfq1}\\
		D^4_{ijrs}\tilde Q^* &= D_l\mathcal{A}_{kr}\left(\frac{dz_l}{dp_s}\right)D_k\mathcal{A}_{ij} + \mathcal{A}_{kr} D^2_{kl}\mathcal{A}_{ij}\left(\frac{dz_l}{dp_s}\right) \nonumber \\
				&= (D_l\mathcal{A}_{kr})\mathcal{A}_{ls}(D_k\mathcal{A}_{ij}) + \mathcal{A}_{kr} (D^2_{kl}\mathcal{A}_{ij}) \mathcal{A}_{ls}, \nonumber
	\end{align}	
and thus
	\begin{equation*}
	\begin{split}
		MTW &=  \mathcal{A}^{rk}\left( \mathcal{A}_{kr} (D^2_{kl}\mathcal{A}_{ij}) \mathcal{A}_{ls} + (D_l\mathcal{A}_{kr})\mathcal{A}_{ls}(D_k\mathcal{A}_{ij}) - \mathcal{A}^{pq}\mathcal{A}_{kp} (D_k\mathcal{A}_{ij}) \mathcal{A}_{ls} (D_l\mathcal{A}_{qr}) \right)\mathcal{A}^{sl} \\
			&= D^2_{kl}\mathcal{A}_{ij}.
	\end{split}
	\end{equation*}		
Therefore, the MTW condition in our case is that:
	\begin{equation}\label{MTW}
		D^2_{z_kz_l}\mathcal{A}_{ij}(z) \xi_i\xi_j\eta_k\eta_l \geq c_0|\xi|^2|\eta|^2
	\end{equation}
for any $\xi, \eta\in\mathbb{R}^d$ with $\xi\bot\eta$, where $c_0>0$ is a constant. 

Next, we shall {formulate} the condition \eqref{MTW} {in terms of} the original {potential} function $\tilde Q$. 
From \eqref{nota}, $\mathcal{A}^{ij}=D^2_{ij}\tilde Q$. By differentiating $I = \mathcal{A}\mathcal{A}^{-1}$, we have
	\begin{align}
		D_k\mathcal{A}_{ij} &= - \mathcal{A}_{ir}\big(D_k\mathcal{A}^{rs}\big)\mathcal{A}_{sj}, \label{mfq2}\\
		D^2_{kl}\mathcal{A}_{ij} &= - \mathcal{A}_{ir}\big(D^2_{kl}\mathcal{A}^{rs}\big)\mathcal{A}_{sj} + 2\mathcal{A}_{im}\big(D_l\mathcal{A}^{mn}\big)\mathcal{A}_{nr}\big(D_k\mathcal{A}^{rs}\big)\mathcal{A}_{sj} \nonumber\\
			&= - \mathcal{A}_{ir}\big(D^2_{kl}\mathcal{A}^{rs}\big)\mathcal{A}_{sj} + 2\big(D_l\mathcal{A}_{ir}\big)\mathcal{A}^{rs}\big(D_k\mathcal{A}_{sj}\big). \nonumber
	\end{align}

Hence, the left hand side of \eqref{MTW} is
	\begin{equation*}
	\begin{split}
		& D^2_{\eta\eta}\mathcal{A}_{\xi\xi} \\ 
		=& - \big(D^2_{\eta\eta}\tilde Q_{ij}\big)(\mathcal{A}\xi)_i(\mathcal{A}\xi)_j + 2 (\mathcal{A}\xi)_i \big(D_\eta\tilde Q_{ir}\big)\tilde Q^{rs}\big(D_\eta\tilde Q_{sj}\big)(\mathcal{A}\xi)_j.
	\end{split}
	\end{equation*}
Letting $\tilde\xi:=\mathcal{A}\xi$ and recalling that $\tilde Q=\tilde Q(z)$ is a scalar function, one has
	\begin{equation}\label{anaMTW}
	\begin{split}
		D^2_{\eta\eta}\mathcal{A}_{\xi\xi} &=  -  D^2_{\eta\eta}\tilde Q_{\tilde\xi\tilde\xi}  + 2 \big(D_\eta\tilde Q_{\tilde\xi r}\big)\tilde Q^{rs}\big(D_\eta\tilde Q_{s\tilde\xi}\big)  \\
		&= -  D^4_{\tilde\xi\tilde\xi\eta\eta}\tilde Q + 2 \left\langle \nabla \tilde Q_{\tilde\xi\eta}, \nabla \tilde Q_{\tilde\xi\eta} \right\rangle_{\mathcal{A}},
	\end{split}
	\end{equation}
where $\langle \xi, \xi \rangle_{\mathcal{A}}:=\xi_i\xi_j\tilde Q^{ij}$.

Combining \eqref{MTW} and \eqref{anaMTW}, we obtain the following result.
\begin{proposition}
The MTW condition in our case can be expressed directly in the terms of the potential function $\tilde Q$:
\begin{itemize}
\item[\textbf{(H2)}] The {potential} function $\tilde Q$ satisfies that for all $\xi, \eta\in\mathbb{R}^d$ with $ \xi \bot \eta$, 
	\begin{equation}\label{newH2}
		\sum_{i,j,k,l,p,q,r,s}\left(\tilde Q_{ijrs} - 2\tilde Q^{pq} \tilde Q_{ijp} \tilde Q_{qrs} \right)\tilde Q^{rk}\tilde Q^{sl} \xi_k\xi_l\eta_i\eta_j \leq -\delta_0|\xi|^2|\eta|^2,
	\end{equation}
where $\{\tilde Q^{ij}\}$ is the inverse matrix of $\{\tilde Q_{ij}\}$, and $\delta_0$ is a positive constant.
When $\delta_0=0$, we call it \textbf{(H2w)}, a weak version of {(H2)}.
\end{itemize} 
\end{proposition}
Comparing with \eqref{dualMTW}, one can see that \eqref{newH2} is in a similar form in spite of the factor $2$.

\vspace{10pt}
\section{Regularity of the potential}\label{s5}

It is well known that in order to guarantee some regularity for equation \eqref{MA**}, one needs some notion of convexity of domains. 
In optimal transport, it has been proved that if the target domain is not $c$-convex, there exist some smooth densities $\rho_0, \rho_T$ such that the optimal mapping is not even continuous, see \cite[\S7.3]{MTW}. For global regularity, one needs both the initial and the target domains to be uniformly $c$-convex in \cite{TW1}.

In our case, the cost function is $c(x,y)=\tilde Q^*(x+y)$.
Similarly to the $c$-convexity in optimal transportation, we introduce the following $q$-convexity for domains. 
\begin{definition}[$q$-exponential map]\label{qexp}
Assume the {potential} function $\tilde Q$ satisfies conditions (H0) and (H1). 
For $x\in\Omega_0$ we define the \emph{q-exponential map} at $x$, denoted by $\mathfrak{I}_x : \R^n \to \R^n$, 
such that
	\begin{equation*}
		\mathfrak{I}_x(z) = \nabla\tilde Q(z) - x.
	\end{equation*}
\end{definition}
Note that this is reminiscent of the mapping $\mathbf{m_2}$ in \eqref{map2}, namely
for $z\in \Om_{T/2}=\nabla\tilde\phi(\Om_0)$, $\mathfrak{I}_x(z) = y \in \Om_T$.
We rename it in order to follow the lines of optimal transportation. 

\begin{definition}\label{dconv}
The domain $\Omega_T$ is \emph{$q$-convex} with respect to $\Omega_0$ if the pre-image $\mathfrak{I}_x^{-1}(\Om_T)$ is convex for all $x\in\Omega_0$, where $\mathfrak{I}_x$ is the $q$-exponential map in Definition \ref{qexp}.

If the pre-image $\mathfrak{I}_x^{-1}(\Om_T)$ is uniformly convex for all $x\in\Omega_0$, then we call $\Omega_T$ is \emph{uniformly $q$-convex} with respect to $\Omega_0$. 

By duality we can also define the (uniform) $q$-convexity for $\Omega_0$ with respect to $\Omega_T$. 
\end{definition}
 
\begin{remark}
\it{(i)} \emph{Similarly to \cite{TW2}, in the smooth case we have an analytic formulation of the $q$-convexity of $\Om_0$ with respect to $\Om_T$. 
Let $\varphi\in C^2(\overline\Om_0)$ be a defining function of $\Om_0$. That is $\varphi=0$, $|\nabla\varphi|\neq0$ on $\p\Om_0$ and $\varphi<0$ in $\Om_0$. $\Omega_0$ is $q$-convex with respect to $\Omega_T$ if
\beq\label{aqcx}
\left[\varphi_{ij}(x) - \tilde {Q^*}^{kl}D_k\tilde Q^*_{ij}(x+y)\varphi_l(x)\right] \geq 0 \qquad \forall x\in\p\Om_0, y\in\Om_T.
\eeq
If the matrix in \eqref{aqcx} is uniformly positive, $\Omega_0$ is uniformly $q$-convex with respect to $\Omega_T$. Note that this analytic formulation is independent of the choice of $\varphi$, and by exchanging $x$ and $y$, we also have the analytic formulation of the $q$-convexity of $\Om_T$ with respect to $\Om_0$.}

\it{(ii)} \emph{Even though $\Omega_0$ and $\Omega_T$ are uniformly $q$-convex, the intermediate domain $\Omega_{T/2}=\nabla\tilde\phi(\Om_0)$ may not be $q$-convex. See \cite{SW} for the counterexamples in the optimal transportation case. }
\end{remark}

\vspace{5pt}

Upon formulating our reconstruction problem to optimal transportation in \S\ref{s3} and assuming appropriate conditions on the {potential} function $\tilde Q$ and the domains, the regularity of the velocity potential $\tilde\phi$ will follow from the established results in optimal transportation. 
In particular, we have the following results that together compose the proof of Theorem \ref{mainthm}.

\begin{theorem}[from \cite{Liu, Lo09}]
Let $\phi$ be the velocity potential in the reconstruction problem \eqref{defIdisc}.
Assume that the {potential} function $\tilde Q$ satisfies conditions (H0), (H1) and (H2), that $\Omega_T$ is $q$-convex with respect to $\Omega_0$.
Assume that $\rho_T\geq c_0$ on $\Om_T$ for some positive constant $c_0$, $\rho_0\in L^p(\Omega_0)$ for some $p>\frac{d+1}{2}$, and the balance condition \eqref{bala} is satisfied.
Then, we have
	\begin{equation*}
		\|\phi\|_{C^{1,\alpha}(\overline\Omega_0)} \leq C,
	\end{equation*}
for some $\alpha\in(0,1)$, where $C$ is a positive constant. 
When $p=\frac{d+1}{2}$, the velocity potential $\phi$ belongs to $C^1(\overline\Om_0)$.
\end{theorem}

\begin{theorem}[from \cite{TW1}]\label{m2}
If furthermore, $\Om_0, \Om_T$ are $C^4$ smooth and uniformly $q$-convex with respect to each other, $\rho_0\in C^2(\overline\Om_0), \rho_T\in C^2(\overline\Om_T)$, then $\phi\in C^3(\overline\Om_0)$, and higher regularity follows from the theory of linear elliptic equations. Particularly, if $\tilde Q, \Om_0, \Om_T, \rho_0, \rho_T$ are $C^\infty$, then the velocity potential $\phi\in C^\infty(\overline\Om_0)$.
\end{theorem}

Similarly to \cite{TW1}, we are able to reduce the condition (H2) to (H2w) in Theorem \ref{m2} by assuming an additional condition called \emph{$c$-boundedness} on $\Om_0$ in \cite{TW1}. Namely, there exists a global barrier function $h$ on $\Om_0$ such that
	\begin{equation*}
		\left[D_{ij}h - \tilde Q^{*lk}\tilde Q^*_{ijl}D_kh\right]\xi_i\xi_j \ge \delta_1|\xi|^2,
	\end{equation*}
for some constant $\delta_1>0$. We refer the reader to \cite{TW1} for more details.

%%%%%%%%%%%%%%%%%%%%%%%%%%%%%%%%%%%
%%%%%%%%%%%%%%%%%%%%%%%%%%%%%%%%%%%
\vspace{10pt}
\section{The mean-field case}\label{s6}

\subsection{Convexity of the space of MTW potentials}
The {potential} function $\tilde Q$ in our two-steps gravitational transport problem is a scalar function defined on the intermediate domain $\Om_{T/2}=\nabla\tilde\phi(\Om_0)$, which only takes effect at time $t=T/2$. 
In the general reconstruction problem considered by the second author in \cite{Lo06}, the gravitational function $p$ actually solves the Poisson equation in the system \eqref{EPs}, and takes effect for all $t\in(0,T)$.  
In fact, at each $t\in(0,T)$, $p$ is a convolution of the density $\rho$ and the Coulomb kernel \eqref{Coulomb}. 
{One can see the convolution as a continuous convex combination with weights $\rho(t,x)$, and a Kernel satisfying (H2) would lead, by convolution, to a potential also satisfying (H2).}
Therefore, in order to study the general case, a natural question one may ask is that: 
\begin{itemize}
\item[ ] \emph{Is the set of {potential} functions satisfying (H2) convex?}
\end{itemize}
The answer is ``No" in general as shown in the following examples.
\begin{example}\label{Ex}
Let the dimension $d=3$ and
	\begin{eqnarray*}
		&& \tilde Q(z) := z_2^2z_3^2 + z_1z_3^2 + z_1z_2^2 + A|z|^2; \\
		&& \tilde Q'(z):= z_2^2z_3^2 - z_1z_3^2 - z_1z_2^2 + A|z|^2,
	\end{eqnarray*}
where $z\in\Om_{T/2}$ that is a bounded domain due to Lemma \ref{C1lem}, and $A$ is a positive constant.

By differentiating we have
	\begin{equation*}
		D^2\tilde Q\ \ (\mbox{or } D^2\tilde Q') =\left[\begin{array}{ccc}
			2A, & \pm 2z_2, & \pm 2z_3 \\
			\pm 2z_2, & 2A + 2z_3^2 \pm 2z_1, & 4z_2z_3 \\
			\pm 2z_3, & 4z_2z_3, & 2A + 2z_2^2 \pm 2z_1
		\end{array}\right].
	\end{equation*}
Choosing $A$ sufficiently large such that $A \gg \diam(\Om_{T/2})$, we have 
	\[ D^2\tilde Q\ \ (\mbox{and } D^2\tilde Q') \ge AI,  \]
where $I$ is the $3\times 3$ identity matrix. 
Hence, $\tilde Q$ and $\tilde Q'$ satisfy conditions (H0) and (H1). 

By further computations one can easily verify that both $\tilde Q$ and $\tilde Q'$ satisfy condition (H2w) as well. 
Choose $\xi=(0,1,0)$, $\eta=(0,0,1)$ and let $MTW(\tilde Q)$ denote the left hand side of \eqref{newH2} for this choice of $\xi,\eta$. 
We have
	\begin{equation*}
		MTW(\tilde Q)\ \ (\mbox{and } MTW(\tilde Q')) \approx \left(D^4_{2233}\tilde Q - 2D^3_{331}\tilde QD^3_{122}\tilde Q\right) \le -2.
	\end{equation*}

However, let $\tilde Q'':=\frac12(\tilde Q+\tilde Q')=z_2^2z_3^2+A|z|^2$. As the third-order terms vanish, one can see that
	\begin{equation*}
		MTW(\tilde Q'') \approx D^4_{2233}\tilde Q'' = 4,
	\end{equation*}
which contradicts with \eqref{newH2}, namely $\tilde Q''$ doesn't satisfy condition (H2).
\end{example}

One can further modify the above example to show that $\tilde Q$, $\tilde Q'$ satisfy (H2) but their convex combination doesn't. 
Define
	\begin{eqnarray*}
		&& \tilde Q(z) := F(z) + T(z) + A|z|^2; \\
		&& \tilde Q'(z):= F(z) - T(z) + A|z|^2,
	\end{eqnarray*}
where $F(z)$ is the fourth-order terms satisfying $F_{ijkl}\xi_i\xi_j\eta_k\eta_l > 0$, and $T(z)$ is the third-order terms such that the quadratic product of $T_{ijs}T_{skl}\xi_i\xi_j\eta_k\eta_l$ is strictly larger than $F_{ijkl}\xi_i\xi_j\eta_k\eta_l$. Then by choosing $A>0$ sufficiently large one has $\tilde Q, \tilde Q'$ satisfy (H0), (H1) and (H2). 
On the other hand, it's easily seen that $\tilde Q'':=\frac12(\tilde Q+\tilde Q')=F(z)+A|z|^2$ doesn't satisfy (H2).

\vspace{5pt}
Inspired by the above example, we may consider the following variety of (H2) by removing the third-order terms and the orthogonality.
\begin{itemize}
\item[\textbf{(H2c)}] The {potential} function $\tilde Q$ satisfies that for all $\xi, \eta\in\mathbb{R}^d$, 
	\begin{equation}\label{H2c}
		\sum_{i,j,k,l,p,q,r,s} \left(D^4_{ijrs}\tilde Q\right) \tilde Q^{rk}\tilde Q^{sl} \xi_k\xi_l\eta_i\eta_j \leq 0.
	\end{equation}
where $\{\tilde Q^{ij}\}$ is the inverse of $\{\tilde Q_{ij}\}$.
\end{itemize}
Notice that the second term on the left hand side of \eqref{newH2} is always non-positive.
Obviously, {(H2c)} implies {(H2w)}, but not the other way around.

\begin{lemma}\label{leconv}
The set $\mathcal{S}$ of {potential} functions $\tilde Q$ satisfying {(H0), (H1), and (H2c)} is convex. 
\end{lemma}

\begin{proof}
For any $\tilde Q\in\mathcal{S}$, define the vector $\tilde\xi$ by $\tilde\xi_i=\tilde Q^{ik}\xi_k$. From (H0) and (H1), one has $|\tilde\xi|\approx|\xi|$, namely there exist two universal constants $C_1, C_2$ such that
	\[ C_1|\xi| \le |\tilde\xi| \le C_2|\xi|. \]
Hence, \eqref{H2c} is equivalent to 
	\begin{equation}\label{H2cc} 
		\sum_{i,j,k,l,p,q,r,s} \left(D^4_{ijrs}\tilde Q\right) \xi_k\xi_l\eta_i\eta_j \leq 0.
	\end{equation}
	
Now the convexity of $\mathcal{S}$ follows naturally from the linearity of \eqref{H2cc}.
In fact, let $\tilde Q, \tilde Q' \in\mathcal{S}$ satisfy \eqref{H2cc}. One can easily check that $\tau\tilde Q+(1-\tau)\tilde Q'$ also satisfies the same inequality \eqref{H2cc}, for all $\tau\in[0,1]$. Therefore, $\tau\tilde Q+(1-\tau)\tilde Q' \in \mathcal{S}$.
\end{proof}

\vspace{5pt}
\subsection{Proof of Theorem \ref{theo:2.5}}

Going back to the notations of \cite{Lo06}, define $J=\rho v$ and then the functional we want to minimize is
\ban
I(\rho, J) = \int_0^T\int_{\R^d}\frac{|J|^2}{2\rho}dxdt + \demi\int_{\R^d}\int_{\R^d}\rho(T/2,x)\kappa(x-y) \rho(T/2,y)dxdy.
\nan
From the assumptions on $\kappa$, one can see that $I(\rho, J)$ is convex in $\rho, J$, and the existence of a unique minimizer follows \cite{Lo06}. It also falls easily that trajectories of minimizers will have constant velocity on $[0, T/2^-)$ and on $(T/2^+,T]$, and that $v(T/2^+,x)-v(T/2^-,x) = \nabla Q(x)$ where $Q$ is given by 
\ban
Q(x) = \int_{\R^d}\rho(T/2,y)\kappa(x-y)dy.
\nan
Therefore, once the $\rho(T/2)$ has been found, the problem becomes equivalent to the non-interacting problem, with $Q$ as above.

\subsection{Proof of Theorem \ref{mainthm2}}

In order to obtain the global regularity, we need to show the following $q$-convexity of domains.
\begin{lemma}\label{cxle}
Assume that $\kappa$ is $C^4$ and convex. Assume $\Om_0$ is convex, and let $\varphi\in C^2(\overline\Om_0)$ be a convex defining function of $\Om_0$. Suppose that
\beq\label{cxas}
\left[\varphi_{ij}(x)+\frac{TM}{8}\kappa_{ijk}(z-w)\varphi_k(x)\right] \geq c_0\delta_{ij} \qquad \forall x\in\p\Om_0, \ \ \forall z, w\in\Om_{T/2}, 
\eeq
for a constant $c_0\geq0$,
where $\kappa_{ijk}$ are partial derivatives in $z$, $M$ is the total mass in \eqref{bala}. 
Then,
$\Om_0$ is uniformly $q$-convex with respect to $\Om_T$ if $c_0>0$; and $\Om_0$ is $q$-convex with respect to $\Om_T$ if $c_0=0$. 
\end{lemma}

\begin{proof}
It suffices to verify the inequality \eqref{aqcx}. In this case, the modified potential is
\beq\label{mfpo}
\tilde Q(z)=\frac{T}{2}\int_{\Om_{T/2}}\bar\rho(w)\kappa(z-w)dw + |z|^2,
\eeq
where $\bar\rho(w)=\rho(T/2,w)$. 
From the proof of Lemma \ref{C1lem}, the intermediate domain $\Om_{T/2}$ is bounded. 
From \cite{Lo09}, the density $\bar\rho$ is bounded, namely $0\leq \bar\rho \leq C$ for a constant $C>0$. 
Since $\kappa$ is convex in $z$, the potential $\tilde Q$ is uniformly convex, namely the matrix
\beq\label{newes1}
\tilde Q_{ij}(z) \geq 2\delta_{ij} \qquad \forall z\in\Om_{T/2}.
\eeq

Now, we convert \eqref{aqcx} in terms of $\tilde Q$. Note that since $\tilde Q^*$ is the Legendre transform of $\tilde Q$, from \eqref{nota}, $\tilde {Q^*}^{kl}(p)=\tilde Q_{kl}(z)=\mathcal{A}^{kl}(z)$, where $p=x+y$. 
Combining \eqref{mfq1} and \eqref{mfq2} one has
\begin{align*}
D^3_{ijk}\tilde Q^*(z) &= \mathcal{A}_{kr}D_r\mathcal{A}_{ij}(z) \\
	&= - \mathcal{A}_{kr} \mathcal{A}_{im}D_r\mathcal{A}^{mn}\mathcal{A}_{nj}(z).
\end{align*}
Hence, \eqref{aqcx} is equivalent to that for any vector $\xi\in \mathbb{S}^{n-1}$,
\beqs
\left[ \varphi_{ij}(x) + \tilde Q^{ir}(D_k\tilde Q_{rs})\tilde Q^{sj}(z)\varphi_k(x)  \right]\xi_i\xi_j \geq 0\qquad\forall x\in\p\Om_0, z\in\Om_{T/2}.
\eeqs
Similarly as for \eqref{anaMTW}, letting $\tilde\xi:=\mathcal{A}\xi$, namely $\tilde\xi_i=\tilde Q^{ir}\xi_r$, one has
\beq\label{newcx}
\left[ \tilde Q_{ir}\varphi_{rs}(x)\tilde Q_{sj} + (D_k\tilde Q_{ij}(z))\varphi_k(x)  \right]\tilde\xi_i\tilde\xi_j \geq 0\qquad\forall x\in\p\Om_0, z\in\Om_{T/2}.
\eeq
Therefore, in order to show $\Om_0$ is $q$-convex with respect to $\Om_T$, it suffices to verify \eqref{newcx}.

By differentiating \eqref{mfpo},
\beqs
D_k\tilde Q_{ij}(z) = \frac{T}{2}\int_{\Om_{T/2}}\bar\rho(w)\kappa_{ijk}(z-w)dw.
\eeqs
From \eqref{bala} and the conservation of mass, one has
\beqs
\int_{\Om_{T/2}}\bar\rho(w)dw = M .
\eeqs
Then, one can easily see that for any $z\in\Om_{T/2}$,
\beq\label{newes2}
\min_{w\in\overline\Om_{T/2}} \frac{TM}{2}\kappa_{ijk}(z-w) \leq D_k\tilde Q_{ij}(z) \leq \max_{w\in\overline\Om_{T/2}} \frac{TM}{2}\kappa_{ijk}(z-w).
\eeq
Since $\varphi$ is convex, combining \eqref{newes1} and \eqref{newes2} into \eqref{newcx} we obtain
\begin{align*}
\left[ \tilde Q_{ir}\varphi_{rs}(x)\tilde Q_{sj} + (D_k\tilde Q_{ij}(z))\varphi_k(x)  \right] &\geq \min_{w\in\overline\Om_{T/2}} 4\left[\varphi_{ij}(x) +\frac{TM}{8}\kappa_{ijk}(z-w)\varphi_k(x) \right] \\
& \geq 4c_0,
\end{align*}
where the last step was due to the assumption \eqref{cxas}.
Note that the matrix $(\tilde Q_{ij}(z))$ is positive definite and bounded for all $z\in\Om_{T/2}$. So, there is a constant $C>0$ such that $C^{-1}|\xi|\leq |\tilde \xi| \leq C|\xi|$ for any vector $\xi\in\R^n$.
Therefore, the lemma is proved.
\end{proof}

\begin{remark}
\it{(i)} \emph{Similarly to the above proof, by exchanging $x$ and $y$, we can also obtain the $q$-convexity of $\Om_T$ with respect to $\Om_0$ under the assumption \eqref{cx2}.}

\it{(ii)} \emph{Note that the convexity assumption on $\Om_0$ can be dropped by imposing a more involved condition on the kernel $\kappa$ comparing to \eqref{cxas}. Here, Lemma \ref{cxle} only provides a sufficient condition for $\Om_0$ to be $q$-convex.}
\end{remark}

\begin{proof}[Proof of Theorem \ref{mainthm2}]
Thanks to Lemma \ref{cxle}, from the assumptions \eqref{cx1}, \eqref{cx2} we know that $\Om_T$ is $q$-convex with respect to $\Om_0$ if $b_1\geq0$, and $\Om_0$, $\Om_T$ are uniformly $q$-convex with respect to each other if $b_0, b_1>0$.
Therefore, it suffices to verify that the modified {potential} function
	\begin{equation*}
		\tilde Q(z) = |z|^2 + \frac{T}{2}\int_{\Omega_{T/2}} \rho(T/2,y) \kappa(z-y) dy
	\end{equation*} 
satisfies the conditions (H0), (H1) and (H2c).
The condition (H0) follows from the smoothness of the kernel $\kappa$;
the condition (H1) follows from the convexity of $\kappa$.
From the assumption that $\kappa$ satisfies the condition (H2c) and Lemma \ref{leconv}, the {potential} function $\tilde Q$ satisfies (H2c).
Hence, the proof follows from Theorem \ref{mainthm}.
\end{proof}

Unfortunately, the Coulomb kernel \eqref{Coulomb} $\kappa(x-y)=\frac{c_d}{|x-y|^{d-2}}$ does not satisfy the condition (H2c). 
By direct computation we can verify that it satisfies the condition (H2). However, since the space of {potential} function satisfying (H2) is not convex, as shown in Example \ref{Ex}, Theorem \ref{mainthm2} does not apply to that case.

\begin{lemma}
The Coulomb kernel $\kappa$ in \eqref{Coulomb} satisfies condition (H2), if $c_d>0$.
\end{lemma}

\begin{proof}
It suffices to consider $\kappa(z)=|z|^{2-d}$.
By differentiation we have
	\beqs
		\nabla\kappa(z) = (2-d)|z|^{-d}z =: p.
	\eeqs
Let $\kappa^*$ be the Legendre transform of $\kappa$, then $\kappa(z)+\kappa^*(p) = z\cdot p$, thus
	\begin{align*}
		\kappa^*(p) &= z\cdot p - \kappa(z) \\
			&= (1-d)|z|^{2-d} \\
			&= \frac{1-d}{(d-2)^{(2-d)/(1-d)}} |p|^{\frac{2-d}{1-d}}.			
	\end{align*}
In the discrete case, from \eqref{cost} we obtain {that} the cost function is
	\beq\label{gravcost}
		c(x,y) = \frac{d-1}{(d-2)^m} |x+y|^{m},\qquad m=\frac{2-d}{1-d}\in(0,1).
	\eeq
By computation in optimal transportation \eqref{cMTW}, we have the LHS of \eqref{MTW} is 
	\begin{align*}
		LHS(z, \xi, \eta) = \frac{d-1}{(d-2)^m}  & \left(\frac{m-2}{m-1}\right) |z|^{-\frac{m}{m-1}}\left\{1-\frac{m}{m-1}\left(\frac{z\cdot\eta}{|z|}\right)^2-m\left(\frac{z\cdot\xi}{|z|}\right)^2\right. \\
			& \left. + \frac{m(3m-2)}{m-1} \left(\frac{z\cdot\xi}{|z|}\right)^2\left(\frac{z\cdot\eta}{|z|}\right)^2 + (m-1)(\xi\cdot\eta)^2 \right\},
	\end{align*}
for any $\xi, \eta$ with $|\xi|=|\eta|=1$. 

Since $0<m<1$, the coefficient $\frac{d-1}{(d-2)^m} \left(\frac{m-2}{m-1}\right) > 0$. For terms in the bracket, we have
	\begin{equation*}
		I := 1 - m\left(\frac{z\cdot\xi}{|z|}\right)^2 + (m-1)(\xi\cdot\eta)^2 \geq 0,
	\end{equation*}
with the equality holds if and only if $\xi=\eta = z/|z|$, and
	\begin{equation*}
	\begin{split}
		I\!I &:= \frac{m}{1-m}\left(\frac{z\cdot\eta}{|z|}\right)^2 - \frac{m(3m-2)}{1-m} \left(\frac{z\cdot\xi}{|z|}\right)^2\left(\frac{z\cdot\eta}{|z|}\right)^2 \\
			&= \frac{m}{1-m}\left(\frac{z\cdot\eta}{|z|}\right)^2\left[ 1 - (3m-2)\left(\frac{z\cdot\xi}{|z|}\right)^2  \right].
	\end{split}
	\end{equation*}

If $m\in(0,\frac23]$, then $(3m-2) \leq 0$, thus
	\begin{equation*}
		I\!I \geq \frac{m}{1-m}\left(\frac{z\cdot\eta}{|z|}\right)^2 \geq 0,
	\end{equation*}
with the last equality holds if and only if $z\perp\eta$. 

If $m\in(\frac23, 1)$, then 
	\beqs
		1 - (3m-2)\left(\frac{z\cdot\xi}{|z|}\right)^2 > 1-(3m-2) = 3-3m>0.
	\eeqs
We also have $I\!I \geq 0$ 	with the last equality holds if and only if $z\perp\eta$. 

Therefore, we obtain
	\begin{equation*}
		LHS(z,\xi,\eta) > 0,
	\end{equation*}
for all unit vectors $\xi, \eta$. 
\end{proof}

%%%%%%%%%%%%%%%%%%%%%%%%
%\newpage


\vspace{10pt}

\begin{thebibliography}{999}

\bibitem{AgCa}
Agueh, M. and Carlier, G.:
\newblock Barycenters in the Wasserstein space.
\newblock {\em SIAM Journal on Mathematical Analysis}, 43 (2011), 904--924.


\bibitem {AGS}
Ambrosio, L.; Gigli, N. and Savar\'e, G.:
\newblock {\em Gradient flows in metric spaces and in the space of probability
  measures}.
\newblock Lectures in Mathematics ETH Z\"urich. Birkh\"auser Verlag, Basel,
  second edition, 2008.

\bibitem {BB} Benamou, J.-D. and Brenier, Y.:
	A computational fluid mechanics solution to the Monge-Kantorovich mass transfer problem.
	\textit{Numer. Math.}, 84 (2000), 375--393.

\bibitem {BCS}
Benamou, J.-D.; Carlier, G. and Santambrogio, F.:
\newblock Variational mean field games.
\newblock In {\em Active Particles, Volume 1}, pages 141--171. Springer, 2017.

	
%\bibitem[Br91]{Br91} Brenier, Y.:
%	Polar factorization and monotone rearrangement of vector-valued functions.
%	\textit{Comm. Pure Appl. Math.}, 44 (1991), 375--417.

\bibitem {Bre}
Brenier, Y.:
\newblock Minimal geodesics on groups of volume-preserving maps and generalized
  solutions of the {E}uler equations.
\newblock {\em Comm. Pure Appl. Math.}, 52 (1999), 411--452.

\bibitem {BF}
Brenier, Y.; Frisch, U.; H\'enon, M.; Loeper, G.; Matarrese, S.; Mohayaee, R. and Sobolevskii, A.:
\newblock Reconstruction of the early universe as a convex optimization
  problem.
\newblock {\em Mon. Not. R. Astron. Soc.}, 346 (2003), 501--524.
	
\bibitem {Caff} Caffarelli, L.:
	Allocation maps with general cost functions, in \textit{Partial Differential Equations and Applications}
	(P. Marcellini, G. Talenti, and E. Vesintini eds). Lecture Notes in Pure and Appl. Math. 177 (1996), 29--35.
	
\bibitem {CW} Chen, S. and Wang, X.-J.:
	Strict convexity and $C^{1,\alpha}$ regularity of potential functions in optimal transportation under condition A3w.
	\textit{J. Differential Equations}, 260 (2016), 1954--1974.
	
\bibitem {FKM} Figalli, A.; Kim, Y.-H. and McCann, R.:
	H\"older continuity and injectivity of optimal maps.
	\textit{Arch. Rational Mech. Anal.}, 209 (2013), 747--795.

\bibitem {FM}
Frisch, U.; Matarrese, S.; Mohayaee, R. and Sobolevski, A.:
\newblock A reconstruction of the initial conditions of the universe by optimal
  mass transportation.
\newblock {\em Nature}, 2002.	
	
	
\bibitem {GM} Gangbo, W. and McCann, R. J.:
	Optimal maps in Monge's mass transport problem.
	\textit{C. R. Acad. Sci. Paris, Series I, Math.}, 321 (1995), 1653--1658.


\bibitem {kanto1}
Kantorovich, L.~V.:
\newblock On mass transportation.
\newblock {\em Dokl. Akad. Nauk. SSSR}, 37 (1942), 227--229.

\bibitem {kanto2}
Kantorovich, L.~V.:
\newblock On a problem of {M}onge.
\newblock {\em Uspekhi Mat.Nauk.}, 3 (1948), 225--226.

\bibitem {KimPass}
Kim, Y.-H. and Pass, B.:
\newblock Wasserstein barycenters over riemannian manifolds.
\newblock {\em Advances in Mathematics}, 307 (2017), 640--683.
	
\bibitem {LM} Lee, P. and McCann R. J.:
	The Ma-Trudinger-Wang curvature for natural mechanical actions.
	\textit{Calc. Var. PDEs}, 41 (2011), 285--299.
	
\bibitem {Liu} Liu, J.:
	H\"older regularity of optimal mappings in optimal transportation.
	\textit{Calc. Var. PDEs}, 34 (2009), 435--451.

\bibitem {Lo06} Loeper, G.:
	The reconstruction problem for the Euler-Poisson system in cosmology.
	\textit{Arch. Rational Mech. Anal.}, 179 (2006), 153--216.
	
\bibitem {Lo09} Loeper, G.:
	On the regularity of solutions of optimal transportation problems.
	\textit{Acta Math.}, 202 (2009), 241--283.


\bibitem {MTW} Ma, X.-N.; Trudinger, N. S. and Wang, X.-J.:
	Regularity of potential functions of the optimal transportation problems.
	\textit{Arch. Ration. Mech. Anal.}, 177 (2005), 151--183.

%\bibitem[McC]{McC} McCann, R. J.:
%	A convexity principle for interacting gases.
%	\textit{Adv. Math.}, 128 (1997), 153--179.
	
\bibitem {MTF} Mohayaee, R.; Tully, R.B. and Frisch, U.:
	Reconstruction of large-scale peculiar velocity fields.
	\textit{Current Issues in Cosmology}, Cambridge University Press, 2006, pp123. 
	
\bibitem {Monge}
Monge, G.:
\newblock M\'emoire sur la theorie des deblais et des remblais.
\newblock In {\em Histoire de l'Acad\'emie Royale des Sciences}. Paris, 1781.
	
\bibitem {R}
Rockafellar, R.~T.:
\newblock {\em Convex analysis}.
\newblock Princeton Landmarks in Mathematics. Princeton University Press,
  Princeton, NJ, 1997.
\newblock Reprint of the 1970 original, Princeton Paperbacks.	
	
\bibitem {SW} Santambrogio, F. and Wang, X.-J.:
	Convexity of the support of the displacement interpolation: Counterexamples.
	\textit{Applied Math. Letters}, 58 (2016), 152--158.
	
\bibitem {TW1} Trudinger, N. S. and Wang, X.-J.:
	On the second boundary value problem for Monge-Amp\`ere type equations and optimal transportation.
	\textit{Ann. Scuola Norm. Sup. Pisa. Cl. Sci.}, 8 (2009), 143--174.

\bibitem {TW2} Trudinger, N. S. and Wang, X.-J.:
	On strict convexity and continuous differentiability of potential functions in optimal transportation.
	\textit{Arch. Rational Mech. Anal.}, 192 (2009), 403--418.

\bibitem {Vi03} Villani, C.:
	\textit{Topics in Optimal Transportation}.
	Graduate Studies in Mathematics, 58.
	Amer. Math. Soc., Providence, RI. 2003.
	
\bibitem {Vi09} Villani, C.:
	\textit{Optimal Transport. Old and New.}
	Grundlehren der Mathematischen Wissenschaften, 338.
	Springer, Berlin-Heidelberg, 2009. 


\end{thebibliography}
\end{document}